\documentclass[preprint,12pt]{elsarticle}

\usepackage{amssymb}

\usepackage{amsthm}
\usepackage{amsmath}
\usepackage{graphicx}
\usepackage{url}
\usepackage{mathrsfs}

\usepackage{newtxtext,newtxmath}

\usepackage{subfigure}
\usepackage{latexsym,amsmath}
\usepackage{euscript}
\usepackage{amsfonts}
\usepackage[noend]{algorithmic}
\usepackage{moreverb}
\usepackage{xcolor}
\usepackage{color}

\usepackage[latin1]{inputenc}
\usepackage[OT2,T1]{fontenc}
\usepackage{algorithm}
\usepackage[noend]{algorithmic}
\usepackage{epstopdf}
\usepackage{siunitx}
\usepackage{amsfonts}
\usepackage{parskip}
\newtheorem{proposition}{Proposition}

\theoremstyle{definition}

\newtheorem{assumption}{Assumption}

\newcommand{\gvn}{\,|\,}

\newcommand{\bb}[1]{\mathbb{#1}}
\renewcommand{\v}[1]{\boldsymbol{#1}}
\newcommand{\m}[1]{\mathrm{#1}}
\renewcommand{\c}[1]{\mathcal{#1}}

\newcommand{\ind}[1]{\mathbb{I}_{\{#1 \}}}

\DeclareMathOperator{\Prob}{\mathbb{P}}
\DeclareMathOperator{\Exp}{\mathbb{E}}

\DeclareMathOperator{\Var}{\mathbb{V}\mathrm{ar}}

\newcommand{\iidDist}{\overset{\mathrm{iid}}{\sim}}

\newcommand{\e}{\mathrm{e}}
\newcommand{\dd}{\mathop{}\mathopen{}\mathrm{d}}

\newcommand{\sign}{\mathrm{sign}}

\journal{Mathematics and Computers in Simulation}

\begin{document}

\begin{frontmatter}

\title{Monte Carlo Estimation of the Density of the Sum of Dependent Random Variables}

\author{Patrick J. Laub}\address{School of Mathematics and Physics, The University of Queensland \\ Department of Mathematics, Aarhus University}
\author{Robert Salomone}\address{School of Mathematics and Physics, The University of Queensland}
\author{Zdravko I. Botev}\address{School of Mathematics and Statistics, University of New South Wales}

\begin{abstract}

We study an unbiased estimator for the density of a sum of random variables  that are simulated from a computer model.  A numerical study on examples with copula dependence is conducted where the proposed estimator performs favourably in terms of variance compared to other unbiased estimators. We provide applications and extensions to the estimation of marginal densities in Bayesian statistics and to the estimation of the density of sums of random variables under Gaussian copula dependence.
\end{abstract}


\end{frontmatter}

\section{Introduction}

Sums of random variables are fundamental to modeling stochastic phenomena. In finance, risk managers need to predict the distribution of a portfolio's future value which is the sum of multiple assets; similarly, the distribution of the sum of an individual asset's returns over time is needed for valuation of some exotic (e.g.\ Asian) options \cite{mcneil2015quantitative,ruschendorf2013mathematical}. In insurance, the probability of ruin (i.e.\ bankruptcy) is determined by the distribution of aggregate losses (sums of individual claims of random size) \cite{klugman2012loss,asmussen2010ruin}. Lastly, wireless system engineers model total interference in a wireless communications network as the sum of all interfering signals (often lognormally distributed) \cite{fischione2007approximation}.

In this article, we consider estimating the probability density function (pdf) of sums of random variables (rvs).  That is, we wish to estimate the pdf of $S = \sum_{k=1}^n X_k$, where $\v X$ is simulated according to the joint pdf $f_{\v X}$. A major motivation for obtaining accurate pdf estimates of a rv is to produce confidence intervals for quantiles. For example, the US Nuclear Regulatory Commission specifies regulations in terms of the ``95/95'' rule, i.e.\ the upper 95\
\[ \widehat{F}_X(x) = \frac1R \sum_{r=1}^R \ind{ X^{[r]} \le x} \quad \text{for } X^{[1]}, \dots, X^{[R]} \iidDist F_X \,, \]
and then the quantile  $\widehat{q}_\alpha=  \widehat{F}_X^{-1}(\alpha)$. In the obvious notation, we then have the convergence in distribution:
\[
\sqrt{R} (\widehat{q}_\alpha - q_\alpha) \stackrel{\mathcal{D}}{\longrightarrow} {\sf N}\left(0,\alpha(1-\alpha)/[f_X(q_\alpha)]^2\right)
\]
 as
$
R \to \infty,
$
where the limiting variance depends on the unknown density $f_X(q_\alpha)$.  Thus, any confidence intervals for $\widehat{q}_\alpha$ require estimation of the density $f_X(q_\alpha)$, which is a highly nontrivial problem.

In general, the pdf of a sum of rvs is only available via an $n$-dimensional convolution.
The convolution usually cannot be computed analytically (except in some special cases, e.g.,\ iid gammas or normals) or numerically via quadrature (unless $n$ is very small).  Approximations have long been applied to this problem in the iid case for large $n$. These include the \emph{central limit theorems}  \cite{Petrov1975}, \emph{Edgeworth expansions}  \cite{Barndorff1989asymptotic}, and \emph{inversion of integral transforms}.

An Edgeworth expansion is a generalization of the CLT, which constructs a non-normal approximation based on the first $K$ moments (equivalently, cumulants) of $S$ (the first term of the edgeworth expansion is the CLT approximation,   which may not be accurate for small $n$).
Moreover, when the summands of $S$ are dependent, then the moment sequence of $S$ is unknown and needs to be estimated (e.g.\ by Monte Carlo), and small errors in the approximation of the higher moments can lead to large errors in the approximation. Hence, the method is not fully deterministic (as it may first appear), and requires careful calibration of the value $K$ to avoid numerical instabilities.

Another common method is to construct the Laplace transform (or characteristic function) of $S$ and numerically invert it, using a method such as those described in \cite{abate2006unified}. However, when the summands are dependent, the Laplace transform of the sum  is unknown, so one has to first estimate it (e.g.\ by Monte Carlo), and then numerically invert this approximation.  Specialized methods have been developed for certain marginals and dependence structures (for example, the sum of lognormals case is considered by \cite{La15}), but an approach for general distributions is still too difficult.

Finally, Monte Carlo estimators such as Conditional Monte Carlo \cite{asmussen2017conditional} and the Asmussen--Kroese estimator \cite{Asmussen2006} utilize details of $\v{X}$'s distribution to produce unbiased estimates with a dimension--independent rate of convergence of $\mathcal O(1/n)$.

The purpose of this work is to  explore an unbiased Monte Carlo estimator for the problem, with a focus on dependent summands.

The estimator is based on treating the pdf estimation problem as a derivative (of the cumulative distribution function) estimation problem. There are several advantages  to the proposed estimator. First, we show that in certain settings it enjoys smaller variance than those based on the Conditional Monte Carlo approach. Secondly, the estimator only requires  evaluation of   the joint pdf up to a (typically unknown) normalizing constant, a situation similar to the application of Markov chain Monte Carlo. As a result of this, the  sensitivity--based approach is useful in estimating posterior marginal densities in Bayesian inference (Section~\ref{ssec:Marginals}). The source code used in this paper is available online \cite{Code}.

In our notation, we use lowercase boldface letters like $\v{c}$, $\v{x}$, $\v{y}$ for non-random vectors and uppercase boldface letters like $\v{X}$ for random vectors, and $\v{1}$ for the vector of 1's. If $\v{X}$ is of length $n$, we write:  $\v{X}=(X_1, \dots, X_n)^\top$. The inner-product is denoted  $\v{x} \cdot \v{y}$. For a differentiable  function $f: \bb R^n\mapsto \bb R$, we write
\[ \nabla f(\v{z}) = \left. \left( \partial f(\v{x}) / \partial x_1 , \dots, \partial f(\v{x}) / \partial x_n \right)^\top \right\vert_{\v{x} = \v{z}}, \]
and use $\nabla_i f(\v{z})$ to denote the $i$'th component of $\nabla f(\v{z})$.

\section{ Sensitivity Estimator} \label{sec:GeneralSetup}
The estimator is derived from a simple application of Likelihood Ratio method \cite{Glynn1990,Reiman1989}, also known as the Score Function method \cite{Rubinstein1986}), that is typically used for derivative estimation of performance measures in Discrete Event Systems. We thus tackle the pdf estimation problem by viewing it as a special type of  sensitivity analysis.
The basic idea appears in \cite[Chapter VII, Example 5.7]{asmussen2007stochastic}, and our contribution is to weaken a technical condition and  use a control variate to reduce variance. Furthermore, for the Gaussian copula case we introduce an novel conditional Monte Carlo estimator that is infinitely smooth (which may be useful in quasi Monte Carlo \cite{asmussen2007stochastic}).

\begin{assumption} \label{assumption}
The random vector $\v{X}$ has a density $f_{\v{X}}$, each $X_i$ is supported either on the entire real line or a half-real line, the gradient $\nabla f_{\v{X}}$ is a continuous function on the  support of $\v{X}$, and we have the integrability condition
$
\bb E\left|\v X \cdot \nabla \log f_{\v{X}}(\v X)\right|<\infty
$ (here $\v X\sim F_{\v X}$). \hfill $\Diamond$
\end{assumption}
This assumption is slightly weaker than the one in \cite[Prop. 3.5 on page 222]{asmussen2007stochastic}, which requires that
$|\frac{\m d}{\m d s}(f_{\v X}(s \v x) s^n)|$ is uniformly bounded by an $f_{\v X}$-integrable function of $\v x$. The proposed estimator is based on the following  simple formula, proved in the appendix.
\begin{proposition} \label{prop:DensityRepresentation}
For the rv $S = \sum_{i=1}^n X_i = \v{1} \cdot \v{X}$ where $\v{X}$ satisfies Assumption~\ref{assumption},
\begin{align} \label{Pushout}
f_S(s) &= \hphantom{-} \frac{1}{s} \Exp \big\{ \ind{\v{1} \cdot \v{X} \le s} [ \v X \cdot \nabla \log f_{\v{X}}(\v X) + n ] \big\}
\end{align}
for any $s \not= 0$. \hfill $\Diamond$
\end{proposition}

It is straightforward to show that \eqref{Pushout} still holds if the indicators $\ind{\cdot}$ are replaced by $-(1-\ind{\cdot})$. This  suggests the pair of  (unbiased)  estimators ($\v X \sim F_{\v X}$):
\[
\begin{split}
\underbrace{\frac{1}{s} \ind{\v{1} \cdot \v{X} \le s} \big[ \v X \cdot \nabla \log f_{\v{X}}\big(\v X\big) + n \big]}_{\textstyle \widehat{f}_{1}(s)}, \text{ and } \underbrace{-\frac{1}{s} \ind{\v{1} \cdot \v{X} > s} \big[ \v X \cdot \nabla \log f_{\v{X}}\big(\v X\big) + n \big]}_{\textstyle \widehat{f}_{2}(s)}.
\end{split}
\]
We make use of both of these estimators by using one as a base estimator and the difference of the two as a control variate (the difference has a known expectation, namely, zero) \cite{asmussen2007stochastic}. In order to ensure the unbiasedness, we may, for example, obtain the control variate coefficient from a pilot (independent) sample, as explained in Section~\ref{sec:NumericalResults}.

\section{ Conditional Monte Carlo Methods}
In the following Sections~\ref{ssec:CondMC} and \ref{ssec:AK}  we describe the Conditional Monte Carlo approach \cite{asmussen2017conditional},  as well as an extension of the Asmussen--Kroese estimator. We then use these methods as benchmarks to illustrate the performance of the proposed estimator in various settings.
\subsection{Conditional Monte Carlo estimator} \label{ssec:CondMC}

The Conditional Monte Carlo estimator \cite{asmussen2017conditional} takes the form
\[ \widehat{f}_{\mathrm{Cond}}(s) = \frac1n \sum_{i=1}^n f_{X_i | \v{X}_{-i}}( s - S_{-i} ) , \quad \v{X} \sim F_{\v{X}},
 \]
where the notation $\v{X}_{-i}$ denotes the vector $\v{X}$ with the $i$-th component removed and $S_{-i}=\v 1 \cdot \v X_{-i}$. This is particularly simple for the independent case, as $f_{X_i | \v{X}_{-i}} = f_{X_i}$.

We now examine the dependent case where $\v{X}$'s dependence structure is given by an Archimedean copula with generator $\psi$; i.e., the cdf yields
\begin{equation*}
\textstyle
\Prob(X_1\le F_{X_1}^{-1}(u_1), \dots, X_n \le F_{X_n}^{-1}(u_n)) =\phi\big( \sum_{i=1}^n \psi(u_i) \big),\quad \v u\in[0,1]^n,
\end{equation*}
where $\phi\equiv \psi^{-1}$ is the functional inverse of $\psi$.
 The  conditional densities of $\v X$ can be calculated from the formula
($\phi^{(n)}$ denotes $n$-th derivative)
\begin{equation} \label{GenCondMC}
f_{X_i | \v{X}_{-i}}(x_i | \v{x}_{-i}) = f_{X_i}(x_i) \psi^{(1)}( F_{X_i}(x_i)) \frac{\phi^{(n)}( \sum_{j=1}^n \psi(F_{X_j}(x_j)) ) }{\phi^{(n-1)}( \sum_{j\not=i} \psi(F_{X_j}(x_j)))} \,.
\end{equation}

Some Archimedean copulas, such as the Clayton and Gumbel--Hougaard copulas, have what is called a Marshall--Olkin representation. An Archimedean copula is in the Marshall--Olkin representation class if  $\phi(s) = \Exp[\e^{-s Z}]$ for some positive rv $Z$ with cdf $F_Z$. Then an $\v X$ with this dependence structure can be simulated via
\begin{equation} \label{MarshallOlkin} \v{X} = \Big(F_{X_1}^{-1}\Big( \phi\Big(\frac{E_1}{Z}\Big) \Big), \dots, F_{X_n}^{-1}\Big( \phi\Big(\frac{E_n}{Z}\Big) \Big) \Big) , \quad E_i \iidDist \mathsf{Exp}(1),\; Z \sim F_Z \,.
\end{equation}
For this case, Asmussen \cite[Proposition~8.3]{asmussen2017conditional} conditions upon the $Z$ as well as $\v{X}_{-i}$ to obtain what we call the \emph{extended Conditional Monte Carlo estimator}
\begin{align} \label{ExtCondMC}
 \widehat{f}_{\mathrm{ExtCond}}(s) &= \frac{1}{n}\sum_{i=1}^n f_{X_i | \v{X}_{-i}, Z}( s - S_{-i} ),
\end{align}
where $f_{X_i | \v{X}_{-i}, Z}(x_i ) = -z \psi'(F_i( x_i )) f_{X_i}( x_i ) \, \e^{ -z \psi(F_i( x_i )) }$ and  $\v{X}$ is given by \eqref{MarshallOlkin}.

We will use this estimator as a benchmark in our comparisons later on.

\subsection{Asmussen--Kroese estimator} \label{ssec:AK}

The Asmussen--Kroese estimator \cite{Asmussen2006} (typically for tail probabilities) is defined as
\[
\widehat{F}_{\mathrm{AK}}(s) = 1 - \sum_{i=1}^n
\overline{F}_{X_i | \v{X}_{-i}}( \max\{ M_{-i}, s - S_{-i} \} ) \]
where: $M_{-i} = \max\{ X_1,\ldots,X_{i-1},X_{i+1},\ldots,X_n\}$  and  $\overline{F}_{X_i | \v{X}_{-i}}(x)=1-F_{X_i | \v{X}_{-i}}(x)$.
Each $\overline{F}_{X_i | \v{X}_{-i}}( \max\{ M_{-i}, s - S_{-i} \} )=
\overline{F}_{X_i | \v{X}_{-i}}( s - S_{-i})$,
whenever $M_{-i}+S_{-i}<s$. Thus,
we can take the derivative of this  piecewise  estimator to obtain
\[
\widehat{f}_{\mathrm{AK}}(s) = \sum_{i=1}^n
f_{X_i | \v{X}_{-i}}( s - S_{-i} ) \ind{M_{-i} + S_{-i} \le s }, \]
which can be viewed as alternative conditional estimator.
When it is applicable, we use the ``extended'' form of this estimator where $f_{X_i | \v{X}_{-i}}$ is replaced with $f_{X_i | \v{X}_{-i}, Z}$ as in Section~\ref{ssec:CondMC}. Notice that the term  $1/n$ in \eqref{ExtCondMC} does not appear here.   We remark that (to the best of our knowledge) this variant of the AK estimator for estimation of a density has not been previously considered.

\section{Numerical Comparisons} \label{sec:NumericalResults}

In this section,

for various distributions of $\v{X}$ we compare: i) our proposed method, ii) the conditional MC estimator, and iii) the Asmussen--Kroese (AK) estimator.

We conduct 3 experiments, each one depicted on  Figures~\ref{Test:Clayton_Weibull} to \ref{Test:Frank_Lognormal} below.
Each experiment uses $R=10^5$ iid replicates of $\v{X}$ which are common to all estimators (our estimator uses the first 5\

For each experiment we display a subplot of the estimated density function, as well as
the estimated standard deviation and (square root of the) \emph{work-normalized relative variance}:
$
 \text{WNRV}(\widehat{f}(x)) = \text{(CPU\_Time)} \times \Var( \widehat{f}(x) ) / (R[\widehat{f}(x)]^2)
$.  Here, CPU\_Time is the (wall) time taken the by method to produce the estimates for the grid of 50 points.

These examples show sums with dependent summands. When the copula has a Marshall--Olkin representation \eqref{MarshallOlkin} we use it to simulate $\v{X}$ and give results for the extended version \eqref{ExtCondMC} of the conditional MC estimator.

Figure~\ref{Test:Clayton_Weibull} considers the sum of dependent identically-distributed heavy-tailed variables. The estimates plot shows us that the estimators basically agree with each other,  as is to be expected when all methods perform well. In terms of WNRV and standard deviation the sensitivity estimator outperforms the others.

\begin{figure}[H]
	\caption{Sum of $n=10$ $\mathsf{Weibull}(0.3, 1)$ random variables with a $\mathsf{Clayton}(1/5)$ copula.}\vspace{3mm}
    \centering
    \includegraphics[scale=1]{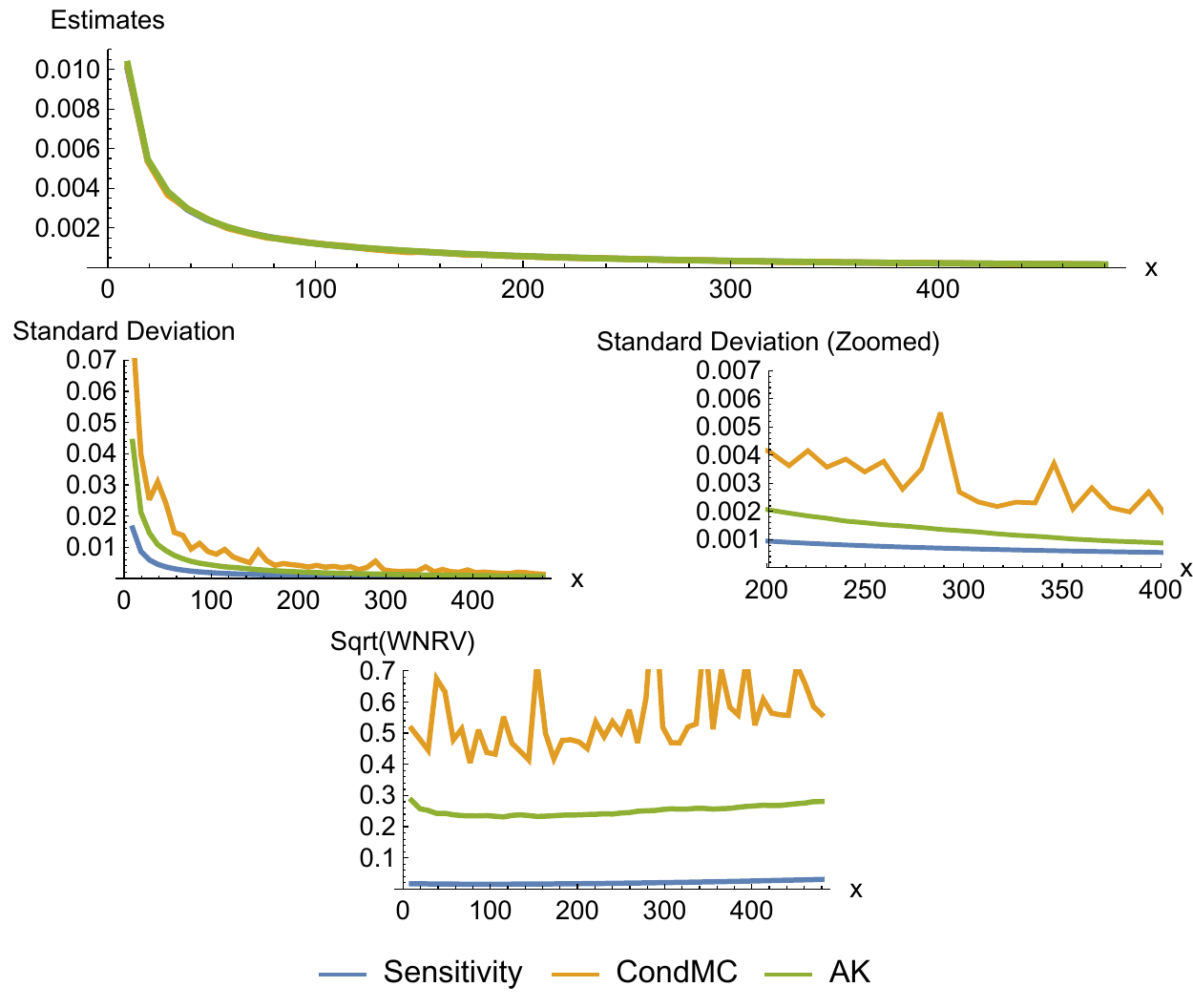}
    \label{Test:Clayton_Weibull}
\end{figure}

Figure~\ref{Test:Gumbel_Exponential} considers a sum of dependent light-tailed variables. The results here are similar to Figure~\ref{Test:Clayton_Weibull}. Again, the sensitivity estimator outperforms the others on WNRV and standard deviation.

\begin{figure}[H]
	\caption{Sum of $n=15$ $\mathsf{Exp}(1)$ random variables with a $\mathsf{GumbelHougaard}(5)$ copula.}\vspace{3mm}
	\centering
	\includegraphics[scale=1]{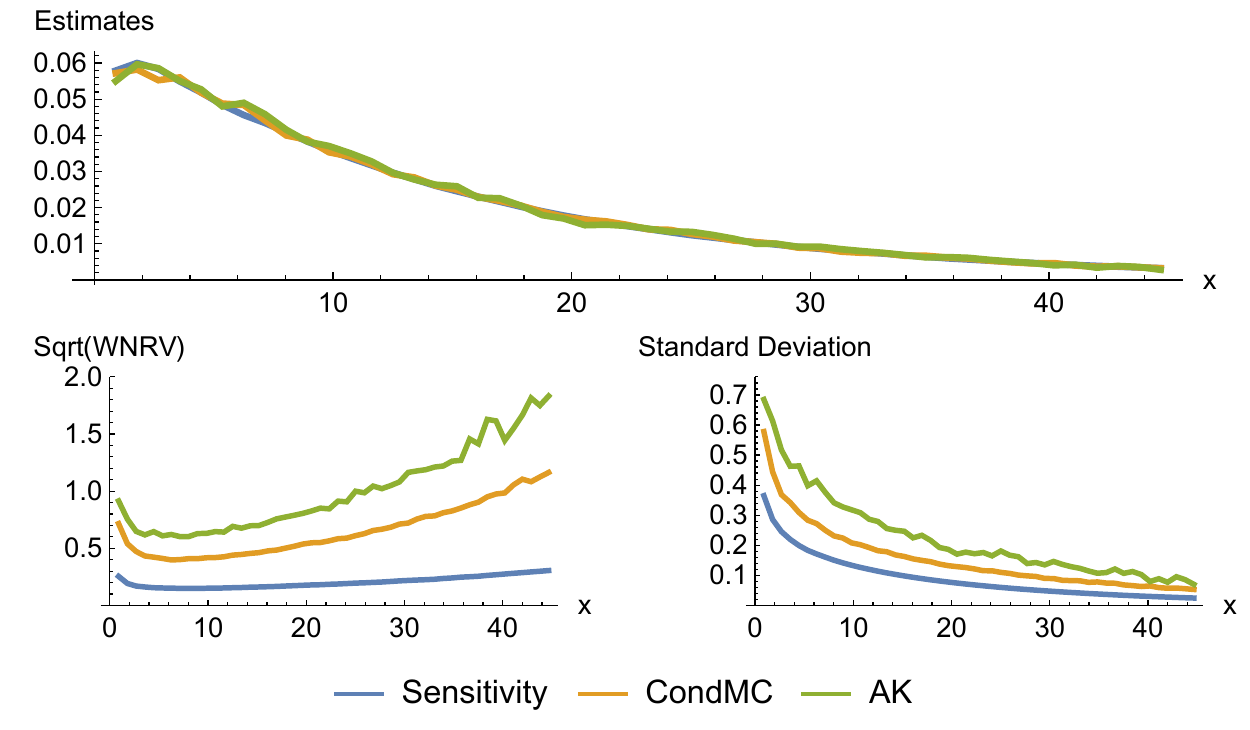}
	\label{Test:Gumbel_Exponential}
\end{figure}

 \begin{figure}[H]
	\caption{Sum of $n=10$ random variables where $X_i \sim \mathsf{Lognormal}(i-10,\sqrt{i})$ with a $\mathsf{Frank}(1/1000)$ copula. The choice of marginals mimic the challenging (and somewhat pathological) example considered in \cite{asmussen2011efficient}.} \vspace{5mm}
    \centering
    \includegraphics[scale=1]{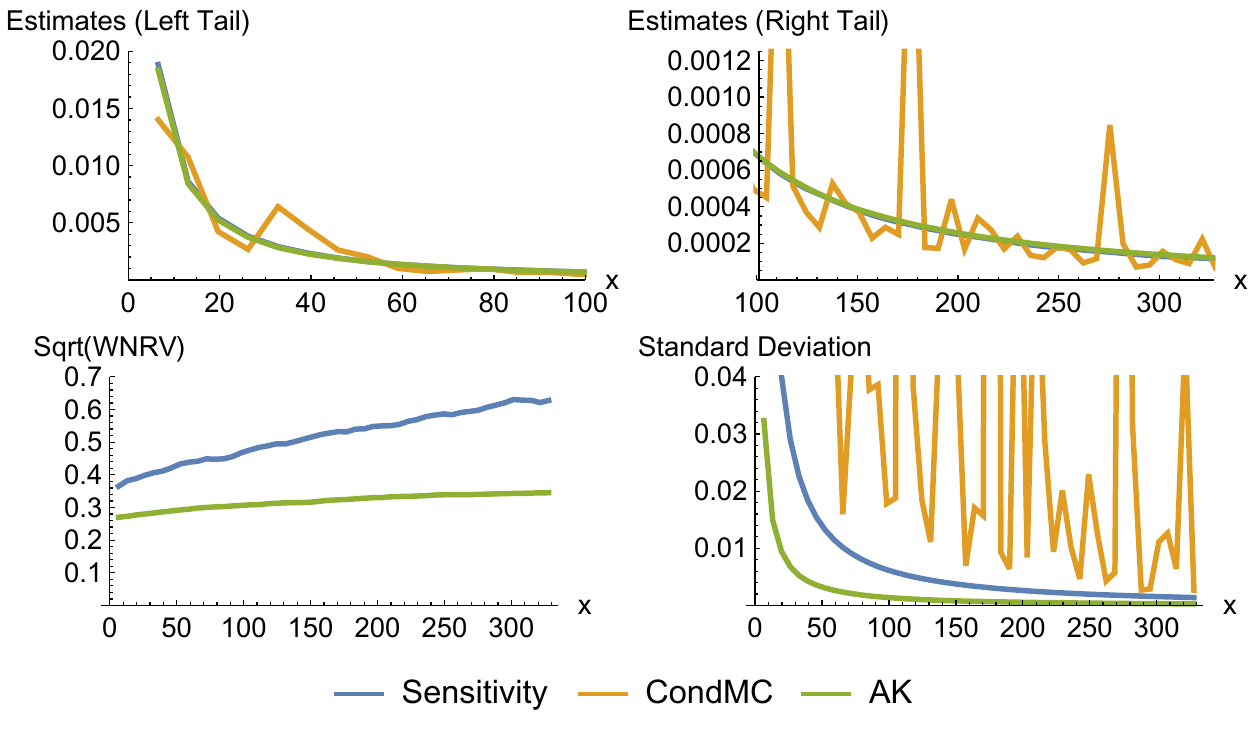}\vspace{3mm}
    \label{Test:Frank_Lognormal}
\end{figure}

Figure~\ref{Test:Frank_Lognormal} shows the sum of dependent heavy-tailed variables. Instead of the standard multivariate lognormal distribution which has a Gaussian copula, we take the Frank copula. The Frank copula is unique among these tests as it is an Archimedean copula which lacks a Marshall--Olkin representation.

\section{Extension to Estimation of Marginal Densities} \label{ssec:Marginals}
One  extension of the sensitivity estimator is in the estimation of marginal densities, which has multiple applications in Bayesian statistics.  For an $\v{X}$ which satisfies Assumption~\ref{assumption}, a similar derivation to the one in Proposition~\ref{prop:DensityRepresentation} gives the following representation of the marginal densities:
    \begin{equation} \label{MarginalExp}
    f_{X_i}(s) = \frac1s \Exp \big\{ \ind{ X_i \le s}  \big(X_i \nabla_i \log f_{\v X}(\v{X}) + 1 \big)  \big\}
    \end{equation}
    for $i=1,\dots,n$, and $s \not=0$.
We use the estimator with associated control variate that is based on \eqref{MarginalExp}. A nice feature of the corresponding estimator is that, due to the presence of the $\nabla \log f_{\v X}(\v x)$ term,  the normalizing constant of $f$ need not be known. As an example, we use Markov Chain Monte Carlo to obtain samples from the posterior density of a Bayesian model, and use these to estimate the posterior marginal pdfs with our  sensitivity estimator.

 We consider the well-known ``Pima Indians'' dataset (standardized), which records a binary response variable  (the incidence of diabetes) for $532$ women, along with seven possible predictors. We specify a Logistic Regression model with predictors: {\em Number of Pregnancies}, {\em Plasma Glucose Concentration}, {\em Body Mass Index}, {\em Diabetes Pedigree Function}, and {\em Age} (see \cite{Friel2012} for justification). The prior is $\v \beta \sim {\sf N}(\v 0, \m I)$, as in \cite{Friel2012}.

To obtain samples from the posterior density, we implement an isotropic Random Walk sampler, using a radially symmetric Gaussian density with $\sigma^2 = 7.5\times 10^{-3}$ (trace plots indicate this choice mixes well for the model).

\begin{figure}
    \centering
    \includegraphics[width=0.8\textwidth]{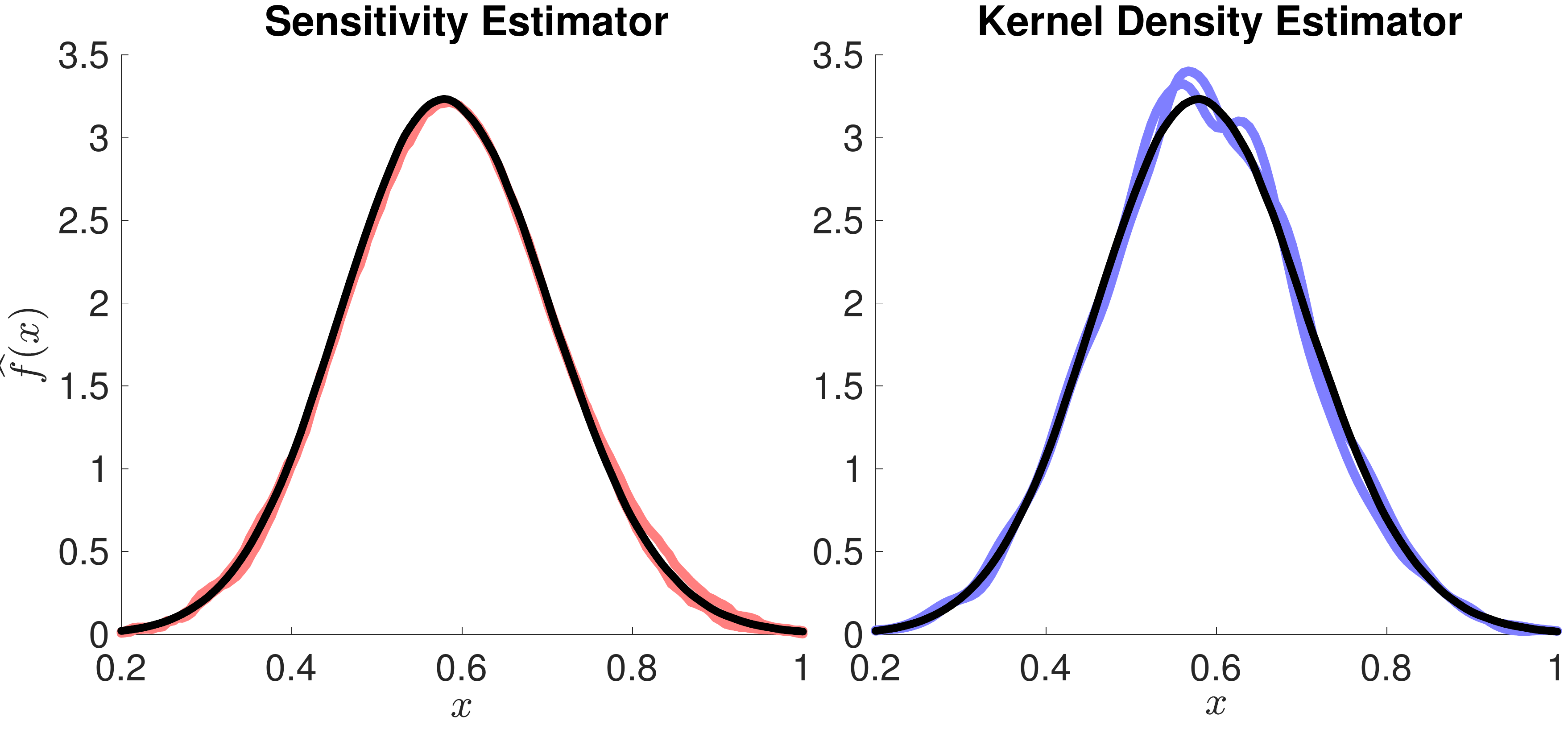}
    \caption{Density estimation of posterior marginal corresponding to the coefficient parameter of the {\em Body Mass Index} predictor variable ( results from two runs are shown).}\label{fig:marginal}
\end{figure}

We ran the Random Walk sampler for $10^3$ steps for burn-in, then used the next $2.5 \times 10^4$ samples (without any thinning) to obtain a KDE, as well as density estimates using our  sensitivity estimator (with control variate).  As a benchmark, we  compare the accuracy with a KDE constructed using every $50$-th  sample from an MCMC chain of length $50\times 5 \times 10^6$. The result of this comparison is depicted in Figure~\ref{fig:marginal}.

As expected, using the same set of samples, the  sensitivity estimator yields a  more accurate estimate than KDE. The reason for the lower accuracy of KDE  in this context is well-known --- a mean square error convergence of
$\c O(n^{-4/5})$, instead of the canonical Monte Carlo rate of $\c O(n^{-1})$,
due to the presence of non-negligible bias in the KDE estimator (see  \cite{botev2010kernel}, for example).

It is important to note that due to the $1/s$ term, the sensitivity estimator can have large variance for very small $s$, even when $F(s)$ or $1 - F(s)$ is not  close to zero. This problem can be resolved with a simple linear shift, as follows.  If one element, say $X_1$, is supported on $\bb{R}$, then $f_S(s) = f_{\widetilde{S}}(s-a)$ for $a \in \mathbb{R}$, where $\widetilde{S} = (X_1+a) + X_2 + \dots + X_n$.  We can then use the original estimator (with shifted values of $s$ and $X_1$) to obtain estimates of the  density of $S$ near or at zero.

\section{Extensions for a Gaussian Copula with Positive Random Variables}
Recall that we wish to estimate the derivative of the cdf
$
F_S(s)=\bb P(\v 1\cdot \v X<s),
$ where $\v X\sim f_{\v X}$. Here we consider random variables that are positive and thus consider $s>0$.   We can then write
$
F_S(s)=\bb P(\v 1\cdot \v Y<1),
$
where $\v Y\sim f_{\v Y}(\v y;s)\equiv f_{\v X}(s\v y)s^n$ (that is, $\v Y$ has the same distribution as $\v X/s$). When $\bb P(\v Y\geq 0)=1$, the nested sequence of events:
\[
\{Y_1\leq 1\}\supseteq\{Y_2\leq 1-Y_1\}
\supseteq \{Y_3\leq 1-Y_1-Y_2\}\supseteq\cdots\supseteq \{Y_n\leq 1-(Y_1+\cdots+Y_{n-1})\}
\]
suggests the possibility of sequentially simulating the components of the vector $\v Y$:
\[
Y_1\rightarrow Y_2\gvn Y_1\rightarrow  Y_3\gvn (Y_1,Y_2) \rightarrow \cdots\rightarrow Y_n\gvn (Y_1,\ldots,Y_{n-1})
\]
In other words,  for $k=1,\ldots,n$ each $Y_k$ is drawn from the conditional density:
\begin{equation}
\label{cond. pdf}
g_k(y_k;s\gvn y_1,\ldots,y_{k-1}):=\frac{f_{\v Y}(y_k;s\gvn y_1,\ldots,y_{k-1})\bb I\{y_k\leq 1-\sum_{j<k}y_j\}}{\alpha_k( y_1,\ldots,y_{k-1};s)},
\end{equation}
where
$
\{\alpha_k\}
$ are normalizing constants that depend on $s$. If we then simulate
\[
\v Y\sim g(\v y;s):=\prod_k g_k(y_k;s\gvn y_1,\ldots,y_{k-1}),
\]
 we obtain the smooth unbiased estimator:
$
 \hat{F}_S(s)=\prod_{k=1}^n \alpha_k(Y_1,\ldots,Y_{k-1};s)
$.
In addition, an application of the likelihood ratio method \cite{Asmussen2006,Glynn1990, Reiman1989,Rubinstein1986} gives us an estimator similar to \eqref{Pushout}, but without the indicator function:
\begin{equation}
\label{smooth est}
\hat{f}_S(s)=\left( \v Y\cdot \nabla \log f_{\v X}(s\v Y)+\frac{n}{s}\right)\prod_{k=1}^n \alpha_k(Y_1,\ldots,Y_{k-1};s)\;.
\end{equation}
Note that the density of  $\v Y$ depends on the specific $s$ at which we estimate $f_S(s)$.  We can use the inverse transform method to remove the dependence on $s$. The inverse transform method allows us to write $\v Y=G(\v U;s)$ for some smooth function $G$ (that depends on $s$) and  a  uniformly distributed vector $\v U\in \bb R^n$ (which does not depend on $s$).

In summary, in order to use estimator \eqref{smooth est},  we must be able to: a) easily simulate from the truncated (or conditional) densities in \eqref{cond. pdf} via the inverse transform method; b) evaluate easily the normalizing constants $\{\alpha_k\}$ of the conditional densities in \eqref{cond. pdf}.

Our finding is that satisfying requirement a) is too difficult for the family of Archimedean copulas. This is because no simple formulas exist for the inverse cdfs of the densities in  \eqref{GenCondMC} and thus numerical root-finding methods are required. Nevertheless,  estimator \eqref{smooth est} is viable for the Gaussian copula under which the vector $[\Phi^{-1}(F_{X_1}(X_1)),\ldots,
\Phi^{-1}(F_{X_n}(X_n))]\sim\mathsf{N}(\v 0,\Sigma)$ for some correlation matrix $\Sigma$. This is because  the conditional densities of the multivariate normal density are also normal and their evaluation requires standard linear algebra manipulations of $\Sigma$.

\begin{figure}[htb]
    \centering
    \includegraphics[width=.7\textwidth]{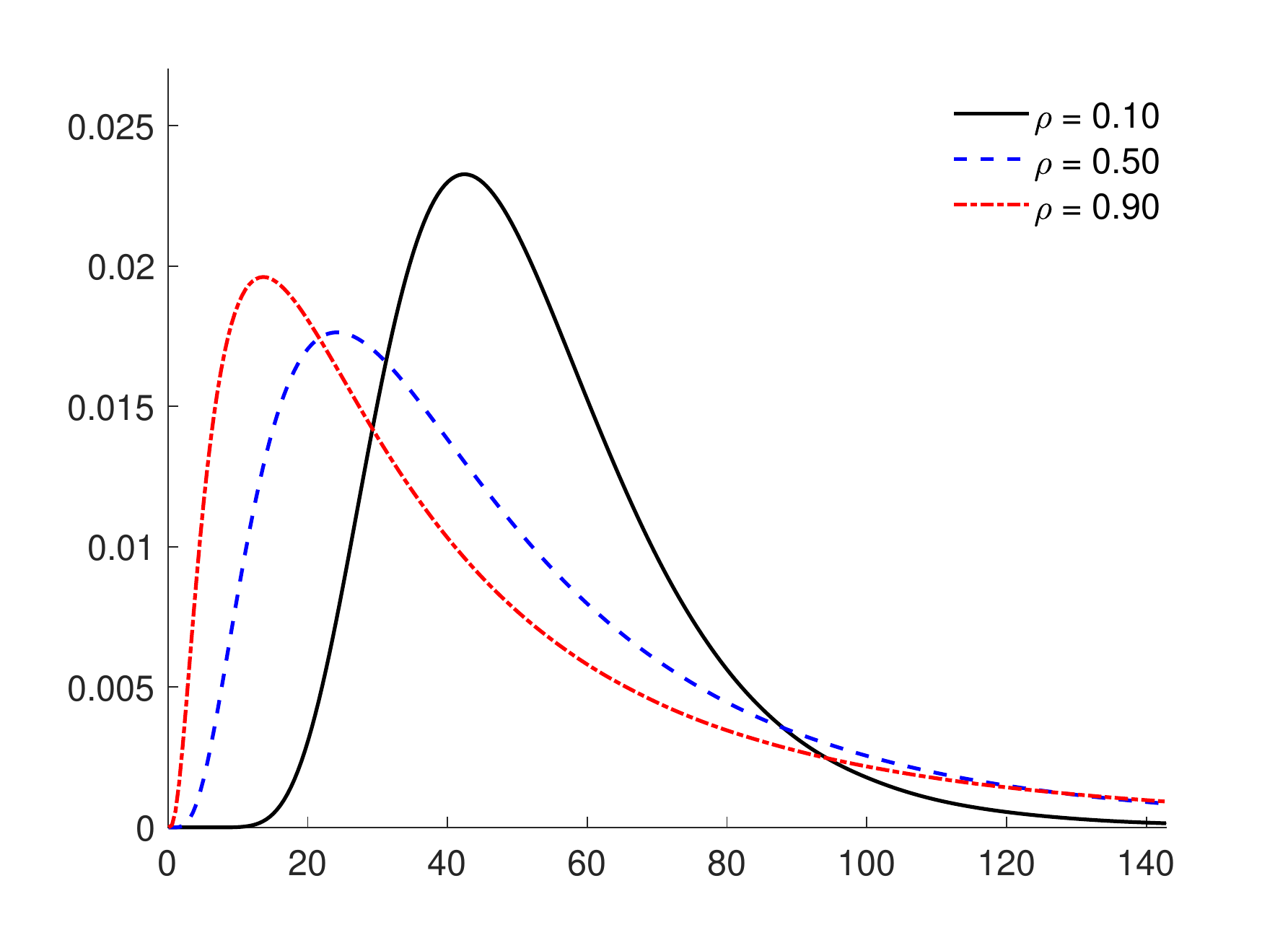}
    \caption{The density of the sum of $32$ dependent lognormal random variables for three correlation values $\rho$. The estimate is based on the average of $R=10^5$ iid replications of \eqref{smooth est}.}\label{fig:lognormal}
\end{figure}

 As an example, Figure~\ref{fig:lognormal} shows the density of the sum of $n=32$ standard lognormal variates
(that is, marginally each $\log(X_k)\sim \mathsf{N}(0,1)$), whose dependence is induced by a Gaussian copula with correlation matrix $\Sigma=\rho\v 1 \v 1^\top+(1-\rho)\m I$ for $\rho\in\{0.1,0.5,0.9\}$.

\section{Conclusion} \label{sec:Conclusion}
In this paper we  derived a sensitivity-based estimator of the pdf of the sum of (dependent) random variables and performed a short numerical comparison. To achieve this we used standard techniques from sensitivity analysis --- we constructed an estimator of the cdf of the sum of random variables and then differentiated this cdf estimator in order to estimate the density. The cdf estimator was constructed using either a change of measure (giving us the likelihood ratio method), or conditional Monte Carlo (as with the Asmussen-Kroese estimator), or both (as  with the Gaussian copula example).

Overall, the numerical comparison indicates that there isn't a single best estimator in all settings. Nevertheless, the proposed sensitivity estimator will likely be preferable in settings where $ \nabla \log f_{\v X}$ can be computed very quickly, and most useful when the conditional Monte Carlo approach is difficult to apply.

As future research we suggest exploring the use of quasirandom numbers in order to further reduce the variance of any smooth density estimators such as \eqref{smooth est}.

\appendix

\section{Proof of Proposition~\ref{prop:DensityRepresentation}}

\footnotesize
There are many ways to derive this formula. One of the simplest is to use the likelihood ratio method (\cite[Ch.VII, (4.1)]{Asmussen2006}, \cite{Glynn1990, Reiman1989,Rubinstein1986}), which requires the interchange of differentiation and integration. A general sufficient condition for this interchange to be valid is given in  \cite[Theorem 1]{LEcuyer1990}.  The proof in this reference uses the dominated convergence theorem, which requires that  $|\frac{\m d}{\m d s} f_{\v{X}}(s\v x)s^n|$ is uniformly bounded by an $f_{\v X}$-integrable function of $\v x$. In our derivation below, we instead use the Fubini-Tonelli theorem, which only requires the integrability of  $|\v x \cdot \nabla \log f_{\v{X}}(\v x)|$ with respect to $f_{\v X}$.

Define the cdf $F_S(s) =\int_{\v 1 \cdot \v x \le s} f_{\v{X}}(\v{x}) \dd \v{x} \,,$
so that the pdf is
$
f_S(s) =\frac{\dd}{\dd s}  F_S(s)
$.
The change of variables  $\v x = s \v y$ yields:
\[
 F_S(s)=\int_{\mathcal{R}_s} f_{\v{X}}(s \v{y}) |s|^n \dd \v{y}\quad s\not =0,
\]
where the notation $\int_{\mathcal{R}_s}$ means $\int_{\v{1}\cdot\v{y}\le 1}$ if $s >0$, else $\int_{\v{1}\cdot\v{y}> 1}$ for $s<0$.

Let $\varphi(s):=\int_{\mathcal{R}_s}\frac{\dd}{\dd s} \left( f_{\v{X}}(s \v{y})  |s|^n\right)\dd \v{y} $. We will use the fact
that $\varphi(s)=f_S(s)$ almost everywhere (i.e.\ except possibly on sets of zero Lebesgue measure) on $s\not\in(-\epsilon,\epsilon)$ for an arbitrarily small $\epsilon>0$.

\enlargethispage{.5cm}
In order to justify  the identity $\varphi(s)=f_S(s)$
(almost everywhere) in the case of
$s>\epsilon$ (similar arguments apply for $s<\epsilon$),
we use the Fubini-Tonelli theorem for exchanging the order of
integration. This exchange
holds under the integrability condition
\begin{equation}
\label{integrability}
\int_\epsilon^s\int_{\v{1}\cdot\v{y}\le 1}\left|\frac{\dd}{\dd t} \left( f_{\v{X}}(t \v{y})  t^n\right)\right|\dd \v{y}\dd t<\infty
\end{equation}
and the existence of a continuous $\nabla f_{\v X}$, both of which follow from  Assumption~\ref{assumption} (verified at the end of this proof). Using the Fubini-Tonelli theorem \cite{rosenthal2006first}
we then write:
\[
\begin{split}
\int_\epsilon^s\varphi(t)\dd t&=\int_\epsilon^s\int_{\v{1}\cdot\v{y}\le 1}\frac{\dd}{\dd t} \left( f_{\v{X}}(t \v{y})  t^n\right)\dd \v{y}\dd t=\int_{\v{1}\cdot\v{y}\le 1}\int_\epsilon^s\frac{\dd}{\dd t} \left( f_{\v{X}}(t \v{y})  t^n\right)\dd t\dd \v{y}\\
&=\int_{\v{1}\cdot\v{y}\le 1} ( f_{\v{X}}(s \v{y})  s^n -
 f_{\v{X}}(\epsilon \v{y})  \epsilon^n)\dd \v{y}=F_S(s)-F_S(\epsilon)
\end{split}
\]
Hence, by the fundamental theorem of Calculus,   $\varphi$ equals the  derivative of $F_S$ up to a set of measure zero.
In other words, $\varphi(s)=f_S(s),s>\epsilon$ almost everywhere. To proceed, we write $\sign(x)= x / |x| = \frac{\dd}{\dd x} |x|$
\begin{align*}
f_S(s)=\varphi(s)
&=\textstyle\int_{\mathcal{R}_s} \left[ \v{y} \cdot \nabla \log f_{\v{X}}(s \v{y}) + \frac{n \, \sign(s)}{|s|} \right] |s|^{n} f_{\v{X}}(s \v y)  \dd \v{y} ,
\end{align*}
so after a change of variables $\v{y} = \v{x} / s$ and using $\sign(x)/|x|=1/x$, we obtain
\begin{align*}
f_S(s)
&= \int_{\v 1 \cdot \v x \le s } \big[ \frac{\v{x}}{s} \cdot \nabla \log f_{\v{X}}(\v{x}) + \frac{n}{s} \big] f_{\v{X}}(\v x)  \dd \v{x}
= \frac{1}{s} \Exp \big\{ \ind{ \v 1 \cdot \v X \le s } [ \v{X} \cdot \nabla \log f_{\v{X}}(\v{X}) + n ] \big\} \,.
\end{align*}
To verify \eqref{integrability}, note that after using the change of variable above, it can be upper bounded by
\[
\textstyle\int_\epsilon^s\frac{1}{t} \Exp \big\{ \ind{ \v 1 \cdot \v X \le t }| \v{X} \cdot \nabla \log f_{\v{X}}(\v{X}) + n |\big\} \dd t\leq (\Exp \big| \v{X} \cdot \nabla \log f_{\v{X}}(\v{X})\big | + n) \int_\epsilon^s\frac{1}{t}\dd t<\infty,
\]
which is bounded by assumption.

\end{document}